\documentclass[11pt]{article} 
\begin{document} 
\newtheorem{prop}{Proposition}[section]
\newtheorem{Def}{Definition}[section] \newtheorem{theorem}{Theorem}[section]
\newtheorem{lemma}{Lemma}[section] \newtheorem{Cor}{Corollary}[section]

\title{\bf Low regularity well - posedness for the one - dimensional Dirac - 
Klein - Gordon system}
\author{{\bf Hartmut Pecher}\\
Fachbereich Mathematik und Naturwissenschaften\\
Bergische Universit\"at Wuppertal\\
Gau{\ss}str.  20\\
D-42097 Wuppertal\\
Germany\\
e-mail {\tt Hartmut.Pecher@math.uni-wuppertal.de}}
\date{}
\maketitle

\begin{abstract}
Local well-posedness for the Dirac -- Klein -- Gordon equations is proven in one 
space dimension, where the Dirac part belongs to $H^{-\frac{1}{4}+\epsilon}$ and 
the Klein - Gordon part to $H^{\frac{1}{4}-\epsilon}$ for $ 0 < \epsilon < 
\frac{1}{4}$ , and global well-posedness, if the Dirac part belongs to the 
charge class $L^2$ and the Klein - Gordon part to $H^k$ with $ 0 < k < 
\frac{1}{2}$ . The proof uses a null structure in both nonlinearities detected 
by d'Ancona, Foschi and Selberg and bilinear estimates in spaces of Bourgain - 
Klainerman - Machedon type.
\end{abstract}

\renewcommand{\thefootnote}{\fnsymbol{footnote}}
\footnotetext{\hspace{-1.8em}{\it 2000 Mathematics Subject Classification:} 
35Q40, 35L70 \\
{\it Key words and phrases:} Dirac -- Klein -- Gordon system,  
well-posedness, Fourier restriction norm method}
\normalsize 
\setcounter{section}{-1}
\section{Introduction}
In this paper we study the Cauchy problem for the Dirac -- Klein -- Gordon 
equations in one space dimension
\begin{eqnarray}
\label{0.1}
-i\beta\frac{\partial}{\partial t} \psi + i \alpha \beta 
\frac{\partial}{\partial 
x} \psi + M \psi & = & g \phi \psi \\
\label{0.2}
\frac{\partial^2}{\partial t^2} \phi -  \frac{\partial^2}{\partial x^2} \phi + 
m^2 \phi & = & \langle \beta \psi,\psi \rangle_{{\bf C}^2}
\end{eqnarray}
with initial data
\begin{equation}
\psi(x,0)  =  \psi_0(x) \,  , \, \phi(x,0)  =  \phi_0(x) \, , \, \frac{\partial 
\phi}{\partial t}(x,0) = \phi_1(x) \, .
\label{0.3}
\end{equation}
Here $\psi$ is a two-spinor field, i.e. $\psi$ has values in ${\bf C}^2$, and 
$\phi$ is a real-valued function. $\alpha$ and $\beta$ are hermitian 
$(2 \times 2)$ -matrices, which fulfill $\alpha^2 = \beta^2 = I$ , $ \alpha 
\beta + \beta 
\alpha = 0$, e.g. we can choose $\alpha = {0\,{-i} \choose i\;\,0}$ , $ \beta = 
{0\,1\choose1\,0}$. $M,m$ and $g$ are real constants.

We are interested in local and global low regularity solutions. The first 
results were obtained by Chadam and Glassey \cite{C},\cite{CG} who proved global 
well-posedness for data $\psi_0 \in H^1$ , $\phi_0 \in H^1$ , $\phi_1 \in  L^2$. 
This result was improved by Bournaveas \cite{B} (cf. also Fang \cite{F}) who 
showed the same results for data $\psi_0 \in L^2$ , $\phi_0 \in H^1,$ $\phi_1 
\in L^2$. Local existence and uniqueness was shown by Fang \cite{F1} for data 
$\psi_0 \in H^{-\frac{1}{4}+\epsilon}, $ $ \phi_0 \in H^{\frac{1}{2} +\delta} $ 
, $ \phi_1 \in H^{-\frac{1}{2} + \delta} $ and $0 < \epsilon \le \frac{1}{4}$ , 
$ 0 < \delta \le 2\epsilon$ . These solutions are global, if $\psi_0 \in L^2$. 
Finally, Bournaveas and Gibbeson \cite{BG} also proved global existence and 
uniqueness for $\psi_0 \in L^2$ , $ \phi_0 \in H^k$ , $ \phi_1 \in H^{k-1}$ for 
$\frac{1}{4} \le k < \frac{1}{2}$. All these global results were obtained by 
using conservation of charge, namely $\int |\psi|^2 \,dx = \int |\psi_0|^2 \, 
dx$ . The energy does not help here because it is not positive definite.

In three space dimensions the best result concerning local well-posedness was 
recently obtained by d'Ancona, Foschi and Selberg \cite{AFS} for data $\psi_0 
\in H^{\epsilon},$ $ \phi_0 \in H^{\frac{1}{2}+\epsilon}$ , $\phi_1 \in 
H^{-\frac{1}{2}+\epsilon}$ with $\epsilon > 0$. This result is arbitrarily close 
to the minimal regularity predicted by scaling ($\epsilon = 0$). Whereas in the 
above mentioned more recent results a null structure of Klainerman - Machedon 
type \cite{KM} for the nonlinearities was already used in one or the other way, 
they showed that the null form $\langle \beta \psi,\psi \rangle$ of the wave 
part is also hidden (by a duality argument) in the Dirac part of the system and 
both nonlinearities can be treated in a similar way. It was also very helpful to 
first diagonalize the system by using the eigenspace projections of the Dirac 
operator (cf. also Beals and B\'{e}zard \cite{BB}). Of course this local result 
does not directly imply a global one.

In the present paper we want to improve the local and global results in one 
space dimension by consequently using this diagonalization of the system and 
applying the Fourier restriction norm method. We are able to show local 
existence and uniqueness for data $\psi_0 \in H^{-l}$, $\phi_0 \in H^k$ , 
$\phi_1 \in H^{k-1}$, provided $l<\frac{1}{4}$, $ k > 0$, $2l+k<1$ , $l+k \le 1$ 
and $k \ge |l|$ . This means that e.g. $k=l=\frac{1}{4}-\epsilon$ is admissible 
as well as $l=0$ , $k=\epsilon$ , thus improving the above mentioned results of 
Fang and Bournaveas -- Gibbeson. These local results easily imply global ones in 
the case $\psi_0 \in L^2$ , $\phi_0 \in H^k$ , $\phi_1 \in H^{k-1}$ for $0 < k < 
\frac{1}{2} $ , using only charge conservation, also improving Bourneveas -- 
Gibbeson.

This paper is organized as follows. First we rewrite the system as a first order 
system in time in diagonal form. We split $\psi$ as the sum $\pi_+(D)\psi + 
\pi_-(D)\psi,$ where $\pi_{\pm}(D)$ are the projections onto the eigenspaces of 
$-i\alpha \frac{\partial}{\partial x}$ , and also split $\phi$ as the sum 
$\phi_+ + \phi_-$ , where the half waves $\phi_+$ and $\phi_-$ are defined in 
the usual way. Then we analyze the components of the nonlinearity $\langle \beta 
\psi,\psi \rangle$ , namely $\langle \beta \pi_{\pm}(D)\psi,\pi_{\pm}(D) \psi' 
\rangle$ for all possible combinations of signs by computing its Fourier symbol. 
It turns out that the symbol is a piecewise constant matrix in Fourier space 
depending only on the signs of the Fourier variables and especially vanishes in 
certain regions. Then we examine which bilinear estimates for the nonlinear 
terms are necessary for local well-posedness in the framework of the $X^{m,b}$ - 
spaces. It turns out that due to duality arguments two similar estimates have to 
be given for $\langle \beta \pi_{\pm}(D)\psi,\pi_{\pm}(D)\psi' \rangle$ in order 
to treat both nonlinearities. These are given in Lemma \ref{Lemma 1} and Lemma 
\ref{Lemma 2}. The local results are summarized in Theorem \ref{Theorem 1}. 
Global existence is a direct consequence of the local results combined with 
charge conservation (Theorem \ref{Theorem 2}).

We construct our solutions in spaces of the type $X^{m,b}_{\varphi}$ defined as 
follows: for an equation of the form $iu_t - \varphi(-i\frac{\partial}{\partial 
x})u = 0$ , where $\varphi$ is a measurable function, let $X^{m,b}_{\varphi}$ be 
the completion of ${\cal S}({\bf R} \times {\bf R})$ with respect to
$$\|f\|_{X^{m,b}_{\varphi}} := \|\langle \xi \rangle^m \langle \tau \rangle^b 
{\cal F}(e^{it\varphi(-i\frac{\partial}{\partial x})} f(x,t))\|_{L^2_{\xi\tau}} 
=\|\langle \xi \rangle^m \langle \tau + \varphi(\xi) \rangle^b 
\tilde{f}(\xi,\tau)\|_{L^2_{\xi\tau}} $$
where $\langle \cdot \rangle := (1+|\cdot|^2)^{\frac{1}{2}}$ , and $\tilde{f}$ 
denotes the Fourier transform of $f$ with respect to $x$ and 
$t$. We also use the time localized spaces $X^{m,b}_{\varphi}[0,T]$ defined by
$$ \|f\|_{X^{m,b}_{\varphi}[0,T]} = \inf_{\tilde{f}_{|[0,T]} = f} 
\|\tilde{f}\|_{X^{m,b}_{\varphi}} \, . $$
The following fact about these spaces is well-known (cf. , e.g. , \cite{GTV}, 
section 2):
if $v$ is a solution of
$$ iv_t-\varphi(-i\frac{\partial}{\partial x})v = F \quad , \quad v(0) = f $$
on a time interval $[0,T]$ , $T \le 1$ , we have for $b'+1 \ge b \ge 0 \ge b' > 
-\frac{1}{2}$ :
\begin{equation}
\label{*'}
\|v\|_{X^{m,b}_{\varphi}[0,T]} \le c \|f\|_{H^m} +cT^{1+b'-b} 
\|F\|_{X^{m,b'}_{\varphi}[0,T]} \, .
\end{equation}
\section{Preliminaries}
First we transform our system (\ref{0.1}),(\ref{0.2}) into a first order system 
(in 
t) in diagonal form.

Multiplying the Dirac equations from the left by $\beta$ leads to
\begin{eqnarray*}
-i\frac{\partial}{\partial t} \psi - i \alpha  \frac{\partial}{\partial x} \psi 
+ 
M \beta \psi & = & g \phi \beta \psi \\
\frac{\partial^2}{\partial t^2} \phi -  \frac{\partial^2}{\partial x^2} \phi + 
m^2 \phi & = & \langle \beta \psi,\psi \rangle_{{\bf C}^2} \, .
\end{eqnarray*}
Following the paper of d'Ancona, Foschi and Selberg we diagonalize the system by 
defining the projections
$$ \pi_{\pm}(\xi) := \frac{1}{2}(I \pm \hat{\xi} \alpha) \, , $$
where $\hat{\xi} := \frac{\xi}{|\xi|}$ . Then we have $\psi = \psi_+ + \psi_-$ 
with $\psi_{\pm} := \pi_{\pm}(D)\psi $ , $ D:= 
\frac{1}{i}\frac{\partial}{\partial x}$ .
Using the identity
$$ -i\alpha \frac{\partial}{\partial x} = \alpha D = |D| \pi_+(D) - |D| \pi_-(D) 
$$
and
\begin{equation}
\label{1}
\pi_{\pm}(\xi)\beta = \frac{1}{2}(I\pm \hat{\xi}\alpha)\beta = \frac{1}{2}(\beta 
\mp \hat{\xi}\beta\alpha) = \beta\pi_{\mp}(\xi)
\end{equation}
we get by application of $\pi_{\pm}(D)$ to the Dirac equation
\begin{eqnarray*}
\pi_{\pm}(D)(-i\frac{\partial}{\partial t}\psi - i\alpha\frac{\partial}{\partial 
x}\psi) & = & \pi_{\pm}(D)(-i\frac{\partial}{\partial t}\psi +|D|\pi_+(D)\psi - 
|D|\pi_-(D)\psi)\\
& = & -i\frac{\partial}{\partial t}\pi_{\pm}(D)\psi \pm |D|\pi_{\pm}(D)\psi \\
& = & -i\frac{\partial}{\partial t}\psi_{\pm} \pm |D| \psi_{\pm}
\end{eqnarray*}
where we also used
$$ \pi_{\pm}(\xi)\pi_{\mp}(\xi) = \frac{1}{4}(I\pm 
\hat{\xi}\alpha)(I\mp\hat{\xi}\alpha)=\frac{1}{4}(I-\hat{\xi}^2 \alpha^2) = 0 $$
and
$$\pi_{\pm}(\xi)\pi_{\pm}(\xi) = 
\frac{1}{4}(I\pm\hat{\xi}\alpha)(I\pm\hat{\xi}\alpha) = \frac{1}{4}(I\pm 
2\hat{\xi}\alpha + \hat{\xi}^2 \alpha^2) = \frac{1}{2}(I\pm \hat{\xi}\alpha) = 
\pi_{\pm}(\xi) $$
(this also implies especially $ \psi_{\pm} = \pi_{\pm}(D) \psi_{\pm} $).

The Dirac equations are thus transformed into
\begin{eqnarray*} 
(-i\frac{\partial}{\partial t} \pm |D|)\psi_{\pm} & = & -M\beta 
\pi_{\mp}(D)(\psi_+ + \psi_-) + g \pi_{\pm}(D)(\phi \beta \psi) \\
& = & -M\beta \psi_{\mp} + g \pi_{\pm}(D)(\phi \beta (\psi_+ + \psi_-)) \, .
\end{eqnarray*}
We also split the function $\phi$ into the sum $\phi = \frac{1}{2}(\phi_+ + 
\phi_-)$ , where
$$ \phi_{\pm} := \phi \pm iA^{-\frac{1}{2}} \frac{\partial \phi}{\partial t} 
\quad , \quad A:= - \frac{\partial^2}{\partial x^2} + m^2 \, $$
Here we assume $m > 0$ and in fact $m=1$. Otherwise we artificially add a term 
$(1-m^2)\phi$ on both sides of the equation at the expense of having an 
additional linear term $c_0 \phi$ in the inhomogeneous part which can easily be 
taken care of. We calculate
\begin{eqnarray*}
(i \frac{\partial}{\partial t} \mp A^{\frac{1}{2}})\phi_{\pm} & = & 
(i\frac{\partial}{\partial t} \mp A^{\frac{1}{2}})(I \pm iA^{-\frac{1}{2}} 
\frac{\partial}{\partial t}) \phi \\
& = & i \frac{\partial}{\partial t} \phi \mp A^{\frac{1}{2}}\phi \mp 
A^{-\frac{1}{2}} \frac{\partial^2}{\partial t^2} \phi - i 
\frac{\partial}{\partial t} \phi \\
& = & \mp A^{-\frac{1}{2}}(A \phi + \frac{\partial^2}{\partial t^2} \phi) \\
& = & \mp A^{-\frac{1}{2}}( \langle \beta \psi,\psi\rangle_{{\bf C}^2} + c_0 
\phi) \, .
\end{eqnarray*}
Thus the Dirac -- Klein -- Gordon system can be rewritten as
\begin{eqnarray}
\label{1.1}
(-i \frac{\partial}{\partial t} \pm |D|) \psi_{\pm} & = & -M\beta \psi_{\mp} + g 
\pi_{\pm}(D)(\frac{1}{2}(\phi_+ + \phi_-)\beta(\psi_+ + \psi_-)) \\
\label{1.2}
(i \frac{\partial}{\partial t} \mp A^{\frac{1}{2}})\phi_{\pm} & = & \mp 
A^{-\frac{1}{2}}(\langle \beta(\psi_+ + \psi_-),\psi_+ + \psi_-\rangle_{{\bf 
C}^2} + c_0(\phi_+ + \phi_-)) \, .
\end{eqnarray}
The initial conditions are transformed into
\begin{equation}
\label{1.3}
\psi_{\pm}(0,x) = \pi_{\pm}(D)\psi_0(x) \quad , \quad \phi_{\pm}(0,x) = 
\phi_0(x) 
\pm iA^{-\frac{1}{2}} \phi_1(x) \, . 
\end{equation}
It turns out that the decisive bilinear form which has to be considered is given 
by $\langle \beta \pi_{[\pm]}(D)\psi,\pi_{\pm}(D)\psi'\rangle_{{\bf C}^2}$ , 
where $[\pm]$ and $\pm$ denote independent signs. We are going to compute its 
symbol. One has to treat
\begin{eqnarray*}
\lefteqn{{\cal F}(\langle \beta\pi_{[\pm]}(D)\psi,\pi_{\pm}(D)\psi'\rangle_{{\bf 
C}^2})(\xi,\tau)} \\ & = & \int \int_* \langle \beta \pi_{[\pm]}(\xi_1) 
\tilde{\psi}(\xi_1,\tau_1),\pi_{\pm}(-\xi_2)\tilde{\psi}'(-\xi_2,-\tau_2)\rangle
_{{\bf C}^2} d\xi_1 d\tau_1 \, , 
\end{eqnarray*}
where * denotes the region $\xi = \xi_1 + \xi_2$ , $\tau = \tau_1 + \tau_2$ .
Because $\pi_{\pm}$ are hermitian and by use of (\ref{1}) and $\pi_+(-\xi) = 
\pi_-(\xi)$ we get
\begin{eqnarray*}
\lefteqn{\langle \beta 
\pi_{[\pm]}(\xi_1)\tilde{\psi}(\xi_1,\tau_1),\pi_{\pm}(-\xi_2)\tilde{\psi}'(-\xi
_2,-\tau_2)\rangle } \\ & = & \langle  \pi_{\pm}(-\xi_2) \beta 
\pi_{[\pm]}(\xi_1)\tilde{\psi}(\xi_1,\tau_1),\tilde{\psi}'(-\xi_2,-\tau_2) 
\rangle \\ 
& = & \langle \beta \pi_{\mp}(-\xi_2) 
\pi_{[\pm]}(\xi_1)\tilde{\psi}(\xi_1,\tau_1),\tilde{\psi}'(-\xi_2,-\tau_2) 
\rangle \\
& = & \langle \beta \pi_{\pm}(\xi_2) 
\pi_{[\pm]}(\xi_1)\tilde{\psi}(\xi_1,\tau_1),\tilde{\psi}'(-\xi_2,-\tau_2) 
\rangle \, .
\end{eqnarray*}
We compute
\begin{eqnarray*}
4 \pi_{\pm}(\xi _2) \pi_+(\xi_1) & = & (I \pm \hat{\xi}_2 \alpha)(I + 
\hat{\xi}_1 \alpha) \\
& = & I \pm \hat{\xi}_1 \hat{\xi}_2\alpha^2 +(\hat{\xi}_1 \pm \hat{\xi}_2)\alpha 
\\
& = & (1\pm \hat{\xi}_1 \hat{\xi}_2)I + (\hat{\xi}_1 \pm \hat{\xi}_2)\alpha \, .
\end{eqnarray*}
If $\xi_1$ and $\xi_2$ have different signs we have $\hat{\xi}_1\hat{\xi}_2 = 
-1$ and $\hat{\xi}_1 = - \hat{\xi}_2$ , thus $\pi_+(\xi_2)\pi_+(\xi_1)$ $= 0$ .
If $\xi_1$ and $\xi_2$ have the same sign we have $\hat{\xi}_1\hat{\xi}_2 = 1$ 
and $\hat{\xi}_1 =  \hat{\xi}_2$, thus $4 \pi_+(\xi_2)\pi_+(\xi_1) = 2(I \pm 
\alpha)$ (+ , if $\xi_1,\xi_2>0$, and --, if $\xi_1,\xi_2<0$).
Similarly $4\pi_-(\xi_2) \pi_+(\xi_1) = 2(I \pm \alpha)$ , if $\xi_1,\xi_2$ have 
different signs, and $\pi_-(\xi_2)\pi_+(\xi_1) = 0$ , if $\xi_1,\xi_2$ have the 
same sign. Thus we have
$$ \langle \beta \pi_{\pm}(\xi_2) 
\pi_{[\pm]}(\xi_1)\tilde{\psi}(\xi_1,\tau_1),\tilde{\psi}'(-\xi_2,-\tau_2) 
\rangle = \langle \gamma 
\tilde{\psi}(\xi_1,\tau_1),\tilde{\psi}'(-\xi_2,-\tau_2)\rangle $$
where
\begin{itemize}
\item in the (+,+) - case and in the (--,--) - case : $\gamma = 
\frac{1}{2}(\beta 
\pm \beta \alpha)$, if  $\xi_1,\xi_2$ have the same sign, and $\gamma = 0$, if 
$\xi_1,\xi_2$ have different signs.
\item in the (+,--) - case and in the (--,+) - case :  $\gamma = 
\frac{1}{2}(\beta \pm \beta \alpha)$, if  $\xi_1,\xi_2$ have different signs, 
and $\gamma = 0$, if $\xi_1,\xi_2$ have the same sign.
\end{itemize}
\section{Local solutions}
We want to construct solutions $\psi_{\pm}$ and $\phi_{\pm}$ in the spaces 
$X^{s,b}_{\pm}$ and $Y^{s,b}_{\pm}$, respectively, defined as follows.
\begin{Def}
$X^{s,b}_{\pm}$ is the completion of ${\cal S}({\bf R}^2)$ with respect to the 
norm
$$\|\psi\|_{X^{s,b}_{\pm}} = \|\langle \xi \rangle^s \langle \tau \pm |\xi| 
\rangle^b \tilde{\psi}(\xi,\tau)\|_{L^2_{\xi\tau}} $$
for ${\bf C}^2$ - valued functions $\psi$. $Y^{s,b}_{\pm}$ is the same space for 
${\bf C}$ - valued functions $\psi$. We also use the localized norms 
$$ \|\psi\|_{X^{s,b}_{\pm}[0,T]} = \inf_{\hat{\psi}_{|[0,t]} = \psi} 
\|\hat{\psi}\|_{X^{s,b}_{\pm}} $$
and similarly $Y^{s,b}_{\pm}[0,T]$ .
\end{Def}
We consider the following (slightly modified) system of integral equations which 
belongs to our Cauchy problem (\ref{1.1}),(\ref{1.2}),(\ref{1.3}).
\begin{eqnarray}
\label{a}
\psi_{\pm}(t) & = & e^{\mp it|D|} \psi_{\pm}(0) \\
\nonumber
& & -i g \int_0^t e^{\mp i(t-s)|D|} \pi_{\pm}(D) (\frac{1}{2} 
(\phi_+(s)+\phi_-(s)) 
\beta(\pi_+(D) \psi_+(s) \\
\nonumber
& & \quad \quad + \pi_-(D) \psi_-(s))) ds + iM \int_0^t e^{\mp i(t-s)|D|} \beta 
\psi_{\mp}(s) ds \\
\label{b}
\phi_{\pm}(t) & = & e^{\mp itA^{\frac{1}{2}}} \phi_{\pm}(0)\\
\nonumber
& & \pm i \int_0^t e^{\mp i(t-s)A^{\frac{1}{2}}} A^{-\frac{1}{2}} \langle 
\beta(\pi_+(D) 
\psi_+(s) + \pi_-(D) \psi_-(s)),\pi_+(D) \psi_+(s) \\
\nonumber
& & \quad \quad + \pi_-(D) \psi_-(s) \rangle ds  \pm i c_0 \int_0^t 
e^{\mp i(t-s)A^{\frac{1}{2}}} A^{-\frac{1}{2}} (\phi_+(s) + \phi_-(s)) ds
\end{eqnarray}
We remark that any solution of this system automatically fulfills 
$\pi_{\pm}(D)\psi_{\pm} = \psi_{\pm}$ , because applying $\pi_{\pm}(D)$ to the 
right hand side of the equations for $\psi_{\pm}$ gives 
$\pi_{\pm}(D)\psi_{\pm}(0) = \pi_{\pm}(D)\pi_{\pm}(D)\psi_0 = \pi_{\pm}(D)\psi_0 
= \psi_{\pm}(0)$ , and the integral terms also remain unchanged, because 
$\pi_{\pm}(D)^2 = \pi_{\pm}(D)$ and $\pi_{\pm}(D)\beta\psi_{\mp}(s) = \beta 
\pi_{\mp}(D)\psi_{\mp}(s) = \beta \psi_{\mp}(s)$ . Thus $\pi_{\pm}(D)\psi_{\pm}$ 
can be replaced by $\psi_{\pm}$ on the right hand sides, thus the system of 
integral equations reduces exactly to the one belonging to our Cauchy problem 
(\ref{1.1}),(\ref{1.2}),(\ref{1.3}).
 
Let now data be given with 
$$ \psi_0 \in H^{-l}({\bf R}) \, , \, \phi_0 \in H^k({\bf R}) \, , \, \phi_1 \in 
H^{k-1}({\bf R}) \, . $$
This implies $ \psi_{\pm}(0) \in H^{-l}({\bf R})$ and $\phi_{\pm}(0) \in 
H^k({\bf R})$ . In order to construct a solution of the integral equations for 
$t \in [0,T]$ with a suitable $T \le 1$ with $ \psi_{\pm} \in X_{\pm}^{-l, 
\frac{1}{2}+\epsilon'}[0,T]$ and $\phi_{\pm} \in 
Y_{\pm}^{k,\frac{1}{2}+\epsilon'}[0,T]$ ($\epsilon' > 0$ small) we only have to 
show the following estimates for the nonlinearities (using standard facts from 
the theory of $X^{s,b}$ - spaces, especially (\ref{*'})).

Concerning (\ref{a}) we need
\begin{equation}
\label{*}
\| \pi_{\pm}(D)(\phi \beta \pi_{[\pm]}(D) \psi) 
\|_{X^{-l,-\frac{1}{2}+2\epsilon'}_{\pm}} \le c \| \phi 
\|_{Y^{k,\frac{1}{2}+\epsilon'}_+} \| \psi 
\|_{X^{-l,\frac{1}{2}+\epsilon'}_{[\pm]}} 
\end{equation}
and the same estimates with  $\| \phi \|_{Y^{k,\frac{1}{2}+\epsilon'}_+}$ 
replaced by $\| \phi \|_{Y^{k,\frac{1}{2}+\epsilon'}_-}$ on the right hand side. 
Again $[\pm]$ denotes a sign independent of $\pm$.

Concerning (\ref{b}) we have to show
\begin{equation}
\label{**}
\| \langle \beta \pi_{[\pm]}(D) \psi, \pi_{\pm}(D) \psi' \rangle 
\|_{Y^{k-1,-\frac{1}{2}+2\epsilon'}_+} \le c 
\|\psi\|_{X^{-l,\frac{1}{2}+\epsilon'}_{[\pm]}} 
\|\psi'\|_{X^{-l,\frac{1}{2}+\epsilon'}_{\pm}}
\end{equation}
and the same estimate with $Y^{k-1,-\frac{1}{2}+2\epsilon'}_+$ replaced by  
$Y^{k-1,-\frac{1}{2}+2\epsilon'}_-$ on the left hand side.

By duality (\ref{*}) is equivalent to
$$ \left| \int \int \langle \pi_{\pm}(D)(\phi \beta \pi_{[\pm]}(D) \psi),\psi' 
\rangle dx dt \right|  \le c \|\phi\|_{Y^{k,\frac{1}{2}+\epsilon'}_+} 
\|\psi\|_{X^{-l,\frac{1}{2}+\epsilon'}_{[\pm]}} 
\|\psi'\|_{X^{l,\frac{1}{2}-2\epsilon'}_{\pm}} \, . $$
The left hand side equals
$$ \left| \int \int \phi \langle \beta \pi_{[\pm]}(D) \psi,\pi_{\pm}(D) \psi' 
\rangle dx dt \right| \, , $$
which can be estimated by
$$  \|\phi\|_{Y^{k,\frac{1}{2}+\epsilon'}_+} \| \langle \beta \pi_{[\pm]}(D) 
\psi , \pi_{\pm} \psi' \rangle \|_{Y^{-k,-\frac{1}{2}-\epsilon'}_+} \, . $$
Thus (\ref{*}) is fulfilled if
\begin{equation}
\label{***}
\| \langle \beta \pi_{[\pm]}(D) \psi , \pi_{\pm}(D) \psi' \rangle 
\|_{Y^{-k,-\frac{1}{2}-\epsilon'}_+} \le c 
\|\psi\|_{X^{-l,\frac{1}{2}+\epsilon'}_{[\pm]}} 
\|\psi'\|_{X^{l,\frac{1}{2}-2\epsilon'}_{\pm}}
\end{equation}
and the same with  $Y^{-k,-\frac{1}{2}-\epsilon'}_+$ replaced by 
$Y^{-k,-\frac{1}{2}-\epsilon'}_-$ on the left hand side.

The linear terms in the integral equations can easily be treated as follows:
\begin{eqnarray*}
\lefteqn{ \|\psi_{\mp}\|_{X^{-l,-\frac{1}{2}+2\epsilon'}_{\mp}[0,T]}} \\
& &  \le \|\psi_{\mp}\|_{L^2([0,T],H^{-l})} \le T^{\frac{1}{2}} 
\|\psi_{\mp}\|_{L^{\infty}([0,T],H^{-l})} \le c T^{\frac{1}{2}} 
\|\psi_{\mp}\|_{X^{-l,\frac{1}{2}+\epsilon'}_{\mp}[0,T]} 
\end{eqnarray*} 
and
\begin{eqnarray*}
\lefteqn{ \|A^{-\frac{1}{2}}\phi_{\pm}\|_{Y^{k,-\frac{1}{2}+2\epsilon'}_{[\pm]} 
[0,T]}} \\
& &  \le \|\phi_{\pm}\|_{L^2([0,T],H^{k-1})} \le 
T^{\frac{1}{2}}\|\phi_{\pm}\|_{L^{\infty}([0,T],H^{k-1})} \le c T^{\frac{1}{2}} 
\|\phi_{\pm}\|_{Y^{k,\frac{1}{2}+\epsilon'}_{\pm}[0,T]} \, . 
\end{eqnarray*}

It remains to prove (\ref{**}) and (\ref{***}).
\begin{lemma}
\label{Lemma 1}
Assume $ l < \frac{1}{4} $ , $ 2l+k < 1 $ and $l+k \le 1$ . Then (\ref{**}) 
holds for a sufficiently small $\epsilon' > 0 $ .
\end{lemma}
{\bf Proof:} We have to show
$$ \left| \int \int \langle \beta \pi_{[\pm]}(D) \psi , \pi_{\pm}(D) \psi' 
\rangle \overline{\phi} dx dt \right| \le c 
\|\phi\|_{Y^{1-k,\frac{1}{2}-2\epsilon'}_+} 
\|\psi\|_{X^{-l,\frac{1}{2}+\epsilon'}_{[\pm]}} 
\|\psi'\|_{X^{-l,\frac{1}{2}+\epsilon'}_{\pm}} \, . $$
The left hand side equals (according to the calculation above)
$$ \left| \int \int_*  \langle \beta \pi_{\pm}(\xi_2) \pi_{[\pm]}(\xi_1) 
\tilde{\psi}(\xi_1,\tau_1) , \tilde{\psi}'(-\xi_2,-\tau_2) \rangle 
\overline{\tilde{\phi}}(\xi,\tau) d\xi_1 d\xi_2 d\tau_1 d\tau_2 \right| \, , $$
where * denotes the region $\xi_1 + \xi_2 = \xi$ , $ \tau_1 + \tau_2 = \tau $ .

Defining now 
\begin{eqnarray*}
\tilde{v}_1(\xi_1,\tau_1) & := & \langle \xi_1 \rangle^{-l} \langle \tau_1 [\pm] 
|\xi_1| \rangle^{\frac{1}{2}+\epsilon'}\tilde{\psi}(\xi_1,\tau_1) \\
\tilde{v}_2(\xi_2,\tau_2) & := & \langle \xi_2 \rangle^{-l} \langle \tau_2 \pm 
|\xi_2| \rangle^{\frac{1}{2}+\epsilon'}\tilde{\psi'}(\xi_2,\tau_2) \\
\tilde{\varphi}(\xi,\tau) & := & \langle \xi \rangle^{1-k} \langle \tau + |\xi| 
\rangle^{\frac{1}{2}-2\epsilon'}\tilde{\phi}(\xi,\tau) 
\end{eqnarray*}
we have
$$ \|\psi\|_{X^{-l,\frac{1}{2}+\epsilon'}_{[\pm]}} = \|v_1\|_{L^2_{xt}} \, , \,
\|\psi'\|_{X^{-l,\frac{1}{2}+\epsilon'}_{\pm}} = \|v_1\|_{L^2_{xt}} \, , \, 
\|\phi\|_{Y^{1-k,\frac{1}{2}-2\epsilon'}_+} = \|\varphi\|_{L^2_{xt}} \, . $$ 
Thus we have to show
\begin{eqnarray*}
\left| \int \int_* \frac{\langle \beta \pi_{\pm}(\xi_2) \pi_{[\pm]}(\xi_1) 
\tilde{v}_1(\xi_1,\tau_1),\tilde{v}_2(-\xi_2,-\tau_2) \rangle \langle \xi_1 
\rangle^l \langle \xi_2 \rangle^l \overline{\tilde{\varphi}}(\xi,\tau)}{ \langle 
\tau_1 [\pm] |\xi_1| \rangle^{\frac{1}{2}+\epsilon'} \langle \tau_2 \mp |\xi_2| 
\rangle^{\frac{1}{2}+\epsilon'} \langle \tau + |\xi| 
\rangle^{\frac{1}{2}-2\epsilon'} \langle \xi \rangle^{1-k}} d\xi_1 d\xi_2 
d\tau_1 d\tau_2 \right| & & \\ 
 \le  c \|v_1\|_{L^2} \|v_2\|_{L^2} \|\varphi\|_{L^2} \, . & & 
\end{eqnarray*}
According to our computations at the end of Section 1 we know: 
in the (+,+) - or (--,--) - case this integral reduces to the region $\xi_1 
\xi_2 > 0$ , whereas in the (+,--) - or (--,+) - case it reduces to $\xi_1 \xi_2 
< 0$ . In any case $\beta \pi_{\pm}(\xi_2) \pi_{[\pm]}(\xi_1)$ is a constant 
matrix in each of the quadrants in the $(\xi_1,\xi_2)$ - plane.

{\bf A.} Let us first consider the (+,--) - or (--,+) - case. Here we have to 
prove
\begin{eqnarray*}
\int \int\limits_{* \atop \xi_1 \xi_2 < 0} 
\frac{|\tilde{v}_1(\xi_1,\tau_1)||\tilde{v}_2(-\xi_2,-\tau_2)| 
|\tilde{\varphi}(\xi,\tau)| \langle \xi_1 \rangle^l \langle \xi_2 \rangle^l}{ 
\langle \tau_1 \pm |\xi_1| \rangle^{\frac{1}{2}+\epsilon'} \langle \tau_2 \pm 
|\xi_2| \rangle^{\frac{1}{2}+\epsilon'} \langle \tau + |\xi| 
\rangle^{\frac{1}{2}-2\epsilon'} \langle \xi \rangle^{1-k}} d\xi_1 d\xi_2 
d\tau_1 d\tau_2 & & \\ 
 \le  c \|v_1\|_{L^2} \|v_2\|_{L^2} \|\varphi\|_{L^2} \, . & & 
\end{eqnarray*}
In this region we have $ |\xi| = ||\xi_1| - |\xi_2|| $ . Define
$$ \sigma_1 = \tau_1 \pm |\xi_1| \, , \, \sigma_2 = \tau_2 \pm |\xi_2| \, , \, 
\sigma = \tau + |\xi| \, . $$
Then we get the decisive algebraic inequality:
\begin{eqnarray}
\lefteqn{2 \min(|\xi_1|,|\xi_2|) \le |\xi_1| + |\xi_2| \mp ||\xi_2| - |\xi_1|| = 
|\xi_1| + |\xi_2| \mp |\xi|}
\label{****} \\ \nonumber
& = & \pm(\tau_1 \pm |\xi_1|) \pm(\tau_2 \pm |\xi_2|) \mp(\tau+|\xi|) = \pm 
\sigma_1 \pm \sigma_2 \mp \sigma \le |\sigma_1| + |\sigma_2| + |\sigma| \, .
\end{eqnarray}
{\bf Case 1:} $|\xi_1| << |\xi_2| $ ($\Rightarrow \, |\xi| \sim |\xi_2|$)\\
(The case $|\xi_2| << |\xi_1| $ can be handled similarly.)\\
We have
$$ \frac{\langle \xi_1 \rangle^l \langle \xi_2 \rangle^l}{\langle \xi 
\rangle^{1-k}} \le c \langle \xi_1 \rangle^l \langle \xi_2 \rangle^{l-1+k} $$
and consider three different cases depending on which of the $\sigma$'s is 
dominant.

{\bf a.} $|\sigma| \ge |\sigma_1| \, , \, |\sigma_2|$ \\
By (\ref{****}) we have $ \langle \sigma \rangle \ge c \langle \xi_1 \rangle$ , 
so that it 
remains to estimate
\begin{eqnarray*}
\int \int_* \frac{|\tilde{v}_1(\xi_1,\tau_1)||\tilde{v}_2(-\xi_2,-\tau_2)| 
|\tilde{\varphi}(\xi,\tau)| \langle \xi_1 \rangle^l \langle \xi_2 
\rangle^{l-1+k}}{ \langle \sigma_1 \rangle^{\frac{1}{2}+\epsilon'} \langle 
\sigma_2 \rangle^{\frac{1}{2}+\epsilon'} \langle \xi_1 
\rangle^{\frac{1}{2}-2\epsilon'}} d\xi_1 d\xi_2 d\tau_1 d\tau_2 & & \\
\le \int \int_* \frac{|\tilde{v}_1(\xi_1,\tau_1)||\tilde{v}_2(-\xi_2,-\tau_2)| 
|\tilde{\varphi}(\xi,\tau)|}{ \langle \sigma_1 \rangle^{\frac{1}{2}+\epsilon'} 
\langle \sigma_2 \rangle^{\frac{1}{2}+\epsilon'} \langle \xi_1 
\rangle^{-2l+\frac{3}{2}-k-2\epsilon'}} d\xi_1 d\xi_2 d\tau_1 d\tau_2 & & \\
= \int \int_* \frac{|\tilde{v}_1(\xi_1,\tau_1)||\tilde{v}_2(-\xi_2,-\tau_2)| 
|\tilde{\varphi}(\xi,\tau)|}{ \langle \sigma_1 \rangle^{\frac{1}{2}+\epsilon'} 
\langle \sigma_2 \rangle^{\frac{1}{2}+\epsilon'} \langle \xi_1 
\rangle^{\frac{1}{2}+\epsilon}} d\xi_1 d\xi_2 d\tau_1 d\tau_2 & &  .
\end{eqnarray*}
Here we used the assumptions $l-1+k \le 0$ , $|\xi_1| \le |\xi_2|$ , $2l+k < 1$. 
$\epsilon' > 0$ is sufficiently small and $\epsilon > 0$.
Forgetting about the factor $\langle \sigma_1 \rangle^{\frac{1}{2}+\epsilon'}$ 
and using Plancherel and H\"older this is bounded by
\begin{eqnarray*}
 \left\| {\cal F}^{-1} \left( \frac{|\tilde{v}_1|}{\langle \xi_1 \rangle^ 
{\frac{1}{2}+\epsilon}} \right) \right\|_{L^2_t L^{\infty}_x}  \left\| {\cal 
F}^{-1} \left( \frac{|\tilde{v}_2(-\xi_2,-\tau_2)|}{\langle \sigma_2 \rangle^ 
{\frac{1}{2}+\epsilon'}} \right) \right\|_{L^{\infty}_t L^2_x} 
\left\|\varphi\right\|_{L^2_tL^2_x} & &   \\
 \le  c \|v_1\|_{L^2_{xt}} \|v_2\|_{L^2_{xt}} \|\varphi\|_{L^2_{xt}} & &
\end{eqnarray*}
by Sobolev's embedding and $X^{0,\frac{1}{2}+\epsilon'}_{\pm} \subset 
L^{\infty}_t L_x^2$ .

{\bf b.} $|\sigma_j|$ ($j=1$ or $j=2$) dominant.\\
This case can be treated similarly by using the estimate
$ \langle \sigma_j \rangle^{\frac{1}{2}+\epsilon'} \ge c  \langle \xi_1 
\rangle^{\frac{1}{2}+\epsilon'} \, . $\\
{\bf Case 2:} $ |\xi_1| \sim |\xi_2| $ .\\
We have
$$ \frac{\langle \xi_1 \rangle^l \langle \xi_2 \rangle^l}{ \langle \xi 
\rangle^{1-k}} \sim \frac{\langle \xi_1 \rangle^{2l}}{\langle \xi 
\rangle^{1-k}} \, .$$

{\bf a.} $ |\sigma| $ dominant.\\
We use (\ref{****}) and get $\langle \sigma \rangle \ge c \langle \xi_1 \rangle$ 
, and moreover , $ l < \frac{1}{4} 
$ , $ \langle \xi_1 \rangle \ge c \langle \xi \rangle $ , and $ 2l+k < 1$ and 
estimate as follows:
\begin{eqnarray*}
\lefteqn{\int \int_* 
\frac{|\tilde{v}_1(\xi_1,\tau_1)||\tilde{v}_2(-\xi_2,-\tau_2)| 
|\tilde{\varphi}(\xi,\tau)| \langle \xi_1 \rangle^{2l-\frac{1}{2}+2\epsilon'}}{ 
\langle \sigma_1 \rangle^{\frac{1}{2}+\epsilon'} \langle \sigma_2 
\rangle^{\frac{1}{2}+\epsilon'} \langle \xi \rangle^{1-k}} d\xi_1 d\xi_2 d\tau_1 
d\tau_2} \\
& \le & \int \int_* 
\frac{|\tilde{v}_1(\xi_1,\tau_1)||\tilde{v}_2(-\xi_2,-\tau_2)| 
|\tilde{\varphi}(\xi,\tau)|}{ \langle \sigma_1 \rangle^{\frac{1}{2}+\epsilon'} 
\langle \sigma_2 \rangle^{\frac{1}{2}+\epsilon'} \langle \xi 
\rangle^{\frac{3}{2}-k-2l-2\epsilon'}} d\xi_1 d\xi_2 d\tau_1 d\tau_2  \\
& = &\int \int_* \frac{|\tilde{v}_1(\xi_1,\tau_1)||\tilde{v}_2(-\xi_2,-\tau_2)| 
|\tilde{\varphi}(\xi,\tau)|}{ \langle \sigma_1 \rangle^{\frac{1}{2}+\epsilon'} 
\langle \sigma_2 \rangle^{\frac{1}{2}+\epsilon'} \langle \xi 
\rangle^{\frac{1}{2}+\epsilon}} d\xi_1 d\xi_2 d\tau_1 d\tau_2  \\
& \le & \|v_1\|_{L^2_{xt}} \left\| {\cal F}^{-1} 
\left(\frac{|\tilde{v}_2(-\xi_2,-\tau_2)|}{\langle \sigma_2 
\rangle^{\frac{1}{2}+\epsilon'}} \right) \right\|_{L^{\infty}_t L^2_x} 
\left\|{\cal F}^{-1} \left(\frac{|\tilde{\varphi}|}{\langle\xi 
\rangle^{\frac{1}{2}+\epsilon}} \right) \right\|_{L^2_t L^{\infty}_x}  \\
& \le & \|v_1\|_{L^2_{xt}}  \|v_2\|_{L^2_{xt}} \|\varphi\|_{L^2_{xt}} 
\end{eqnarray*}

{\bf b.} The cases $|\sigma_1|$ or $|\sigma_2|$ dominant are handled similarly.

{\bf B.} Let us next consider the (+,+) - or (--,--) - case. We have to prove
\begin{eqnarray*}
\int \int\limits_{* \atop \xi_1 \xi_2 > 0} 
\frac{|\tilde{v}_1(\xi_1,\tau_1)||\tilde{v}_2(-\xi_2,-\tau_2)| 
|\tilde{\varphi}(\xi,\tau)| \langle \xi_1 \rangle^l \langle \xi_2 \rangle^l}{ 
\langle \tau_1 \pm |\xi_1| \rangle^{\frac{1}{2}+\epsilon'} \langle \tau_2 \mp 
|\xi_2| \rangle^{\frac{1}{2}+\epsilon'} \langle \tau + |\xi| 
\rangle^{\frac{1}{2}-2\epsilon'} \langle \xi \rangle^{1-k}} d\xi_1 d\xi_2 
d\tau_1 d\tau_2 & & \\ 
 \le  c \|v_1\|_{L^2} \|v_2\|_{L^2} \|\varphi\|_{L^2} \, . & & 
\end{eqnarray*}
In this region we have $|\xi| = |\xi_1| + |\xi_2|$ . Assuming w.l.o.g. $|\xi_2| 
\ge |\xi_1|$ we have $|\xi| \sim |\xi_2|$ and also
$$\frac{ \langle \xi_1 \rangle^l \langle \xi_2 \rangle^l}{\langle \xi 
\rangle^{1-k}} \le \frac{c \langle \xi_1 \rangle^l}{\langle \xi_2 
\rangle^{1-k-l}} \le c \langle \xi_1 \rangle^{2l+k-1} \le c \langle \xi_1 
\rangle^{-\epsilon}$$
by our assumptions $l+k \le 1$ and $2l+k<1$ . Moreover, defining
$$ \sigma_1 = \tau_1 \pm |\xi_1| \, , \, \sigma_2 = \tau_2 \mp |\xi_2| \, , \, 
\sigma = \tau + |\xi| \, , $$
we get
\begin{eqnarray*}
2 \min(|\xi_1|,|\xi_2|) & \le & \mp |\xi_1| \pm |\xi_2| + |\xi_1| + |\xi_2| = 
-(\tau_1 \pm |\xi_1|) - (\tau_2 \mp |\xi_2|) + \tau + |\xi| \\
& = & - \sigma_1 - \sigma_2 + \sigma \le |\sigma_1| + |\sigma_2| + |\sigma| \, .
\end{eqnarray*}

{\bf a.} $|\sigma|$ dominant. \\
This implies $ \langle \sigma \rangle \ge c \langle \xi_1 \rangle $ so that we 
estimate for sufficiently small $\epsilon' > 0$ :
\begin{eqnarray*}
\lefteqn{\int \int_* 
\frac{|\tilde{v}_1(\xi_1,\tau_1)||\tilde{v}_2(-\xi_2,-\tau_2)||\tilde 
\varphi(\xi,\tau)|}{\langle \xi_1 \rangle^{\epsilon} \langle \xi_1 
\rangle^{\frac{1}{2}-2\epsilon'} \langle \sigma_1 
\rangle^{\frac{1}{2}+\epsilon'} \langle \sigma_2 
\rangle^{\frac{1}{2}+\epsilon'}}} \\
& \le & \left\|{\cal F}^{-1} \left(\frac{|\tilde{v}_1|}{\langle \xi_1 
\rangle^{\frac{1}{2}+\epsilon-2\epsilon'}}\right) \right\|_{L^2_t L^{\infty}_x} 
\left\|{\cal F}^{-1} \left(\frac{|\tilde{v}_2(-\xi_2,-\tau_2)|}{\langle \sigma_2 
\rangle^{\frac{1}{2}+\epsilon'}}\right) \right\|_{L^{\infty}_t L^2_x} 
\left\|\varphi\right\|_{L^2_t L^2_x} \\
& \le & c \|v_1\|_{L^2_{xt}} \|v_2\|_{L^2_{xt}} \|\varphi\|_{L^2_{xt}} \, .
\end{eqnarray*}

{\bf b.} The cases $|\sigma_1|$ or $|\sigma_2|$ dominant are handled 
similarly.\\[1em]
{\bf Remark:} The modified estimate (\ref{**}) with 
$Y^{k-1,-\frac{1}{2}+2\epsilon'}_+$ replaced by 
$Y^{k-1,-\frac{1}{2}+2\epsilon'}_-$ is proven in the same way replacing $\sigma 
= \tau + |\xi|$ by $\sigma = \tau - |\xi|$ everywhere. One just has to show that 
the decisive algebraic inequality $ 2 \min(|\xi_1|,|\xi_2|) \le |\sigma_1| + 
|\sigma_2| + |\sigma|$ still holds true. This can easily be seen as follows: in 
Part A of the proof we estimate
\begin{eqnarray*}
\lefteqn{2 \min(|\xi_1|,|\xi_2|) \le |\xi_1| + |\xi_2| \pm \left| |\xi_1| - 
|\xi_2|\right| = |\xi_1| + |\xi_2| \pm |\xi|} \\
& & =  \pm(\tau_1 \pm|\xi_1|) \pm(\tau_2 \pm |\xi_2|)\mp(\tau-|\xi|)
 =  \pm \sigma_1 \pm \sigma_2 \mp \sigma \le |\sigma_1| + |\sigma_2| + |\sigma| 
\, ,
\end{eqnarray*}
and in Part B we get
\begin{eqnarray*}
\lefteqn{2 \min(|\xi_1|,|\xi_2|) \le \pm|\xi_1| \mp |\xi_2| + |\xi_1| + |\xi_2| 
= \pm|\xi_1| \mp |\xi_2| + |\xi|} \\
& & =  (\tau_1 \pm|\xi_1|) +(\tau_2 \mp |\xi_2|) - (\tau-|\xi|)
 =  \sigma_1 + \sigma_2 - \sigma \le |\sigma_1| + |\sigma_2| + |\sigma| \, .
\end{eqnarray*}
\begin{lemma}
\label{Lemma 2}
Assume $ k \ge |l| $ and $ k > 0$ . Then (\ref{***}) holds for a sufficiently 
small $\epsilon' > 0$ .
\end{lemma}
{\bf Proof:} Arguing as in the previous proof we have to show
\begin{eqnarray*}
\left| \int \int_* \frac{\langle \beta \pi_{\pm}(\xi_2) \pi_{[\pm]}(\xi_1) 
\tilde{v}_1(\xi_1,\tau_1),\tilde{v}_2(-\xi_2,-\tau_2) \rangle \langle \xi_1 
\rangle^l \overline{\tilde{\varphi}}(\xi,\tau)}{ \langle 
\tau_1 [\pm] |\xi_1| \rangle^{\frac{1}{2}+\epsilon'} \langle \tau_2 \mp |\xi_2| 
\rangle^{\frac{1}{2}-2\epsilon'} \langle \tau + |\xi| 
\rangle^{\frac{1}{2}+\epsilon'} \langle \xi \rangle^k \langle \xi_2 \rangle^l} 
d\xi_1 d\xi_2 
d\tau_1 d\tau_2 \right| & & \\ 
 \le  c \|v_1\|_{L^2} \|v_2\|_{L^2} \|\varphi\|_{L^2} \, . & & 
\end{eqnarray*}

{\bf A:} Consider first the (+,--) - or (--,+) - case. One has to show
\begin{eqnarray*}
\int \int\limits_{* \atop \xi_1 \xi_2 < 0} 
\frac{|\tilde{v}_1(\xi_1,\tau_1)||\tilde{v}_2(-\xi_2,-\tau_2)| 
|\tilde{\varphi}(\xi,\tau)| \langle \xi_1 \rangle^l}{ 
\langle \tau_1 \pm |\xi_1| \rangle^{\frac{1}{2}+\epsilon'} \langle \tau_2 \pm 
|\xi_2| \rangle^{\frac{1}{2}-2\epsilon'} \langle \tau + |\xi| 
\rangle^{\frac{1}{2}+\epsilon'} \langle \xi \rangle^k \langle \xi_2 \rangle^l} 
d\xi_1 d\xi_2 
d\tau_1 d\tau_2 & & \\ 
 \le  c \|v_1\|_{L^2} \|v_2\|_{L^2} \|\varphi\|_{L^2} \, . & & 
\end{eqnarray*}
In this region we have $ |\xi| = ||\xi_1| - |\xi_2|| $ . Define
$$ \sigma_1 = \tau_1 \pm |\xi_1| \, , \, \sigma_2 = \tau_2 \pm |\xi_2| \, , \, 
\sigma = \tau + |\xi| \, . $$
Again as in the previous proof (cf. (\ref{****})) :
\begin{equation}
\label{*****}
2 \min(|\xi_1|,|\xi_2|) \le |\sigma_1| + |\sigma_2| + |\sigma| \, . 
\end{equation}
{\bf Case 1:} $|\xi_1| << |\xi_2| $ ($\Rightarrow \, |\xi| \sim |\xi_2|$)\\
We have 
$$ \frac{\langle \xi_1 \rangle^l}{\langle \xi \rangle^k \langle \xi_2 \rangle^l} 
\sim \frac{\langle \xi_1 \rangle^l}{\langle \xi_2 \rangle^{k+l}} \, .$$
In the $|\sigma_2|$ - dominant case it remains to estimate, using $\langle 
\xi_1\rangle \le c \langle \sigma_2 \rangle$ , $k+l \ge 0,$ $ k > 0$ and 
$|\xi_2| \ge |\xi_1|$ :
\begin{eqnarray*}
\lefteqn{\int \int_* 
\frac{|\tilde{v}_1(\xi_1,\tau_1)||\tilde{v}_2(-\xi_2,-\tau_2)| 
|\tilde{\varphi}(\xi,\tau)| \langle \xi_1 \rangle^l}{ 
\langle \sigma_1 \rangle^{\frac{1}{2}+\epsilon'}
 \langle \sigma 
\rangle^{\frac{1}{2}+\epsilon'} \langle \xi_1 \rangle^{\frac{1}{2}-2\epsilon'} 
\langle \xi_2 \rangle^{k+l}}  d\xi_1 d\xi_2 d\tau_1 
d\tau_2} \\
& \le & \int \int_* 
\frac{|\tilde{v}_1(\xi_1,\tau_1)||\tilde{v}_2(-\xi_2,-\tau_2)| 
|\tilde{\varphi}(\xi,\tau)|}{ \langle \sigma_1 \rangle^{\frac{1}{2}+\epsilon'} 
\langle \sigma \rangle^{\frac{1}{2}+\epsilon'} \langle \xi_1 
\rangle^{k+\frac{1}{2}-2\epsilon'}} d\xi_1 d\xi_2 d\tau_1 d\tau_2 \\
& \le & \left\| {\cal F}^{-1} \left( \frac{|\tilde{v}_1|}{\langle \xi_1 \rangle^ 
{k+\frac{1}{2}-2\epsilon'}} \right) \right\|_{L^2_t L^{\infty}_x} \|v_2\|_{L^2_t 
L^2_x} \left\| {\cal 
F}^{-1} \left( \frac{|\tilde{\varphi}|}{\langle \sigma \rangle^ 
{\frac{1}{2}+\epsilon'}} \right) \right\|_{L^{\infty}_t L^2_x} 
   \\
& \le & c \|v_1\|_{L^2_{xt}} \|v_2\|_{L^2_{xt}} \|\varphi\|_{L^2_{xt}} \, .   
\end{eqnarray*}
The regions where $|\sigma|$ or $|\sigma_1|$ are dominant are treated 
similarly.\\
{\bf Case 2:} $|\xi_2| << |\xi_1| $ ($\Rightarrow \, |\xi| \sim |\xi_1|$)\\
Using
$$ \frac{\langle \xi_1 \rangle^l}{\langle \xi \rangle^k \langle \xi_2 \rangle^l} 
\sim \frac{1}{\langle \xi_1 \rangle^{k-l}\langle \xi_2 \rangle^l}$$
and (\ref{*****}) we have to estimate in the $|\sigma_2|$ - dominant case, using 
$k-l \ge 0$ , $ k > 0$ and $|\xi_1| \ge |\xi_2|$ :
\begin{eqnarray*}
\lefteqn{\int \int_* 
\frac{|\tilde{v}_1(\xi_1,\tau_1)||\tilde{v}_2(-\xi_2,-\tau_2)| 
|\tilde{\varphi}(\xi,\tau)|}{ 
\langle \sigma_1 \rangle^{\frac{1}{2}+\epsilon'}
 \langle \sigma 
\rangle^{\frac{1}{2}+\epsilon'} \langle \xi_1 \rangle^{k-l} \langle \xi_2 
\rangle^{l+\frac{1}{2}-2\epsilon'}}  d\xi_1 d\xi_2 d\tau_1 
d\tau_2} \\
& \le & \int \int_* 
\frac{|\tilde{v}_1(\xi_1,\tau_1)||\tilde{v}_2(-\xi_2,-\tau_2)| 
|\tilde{\varphi}(\xi,\tau)|}{ 
\langle \sigma_1 \rangle^{\frac{1}{2}+\epsilon'}
 \langle \sigma 
\rangle^{\frac{1}{2}+\epsilon'} \langle \xi_2 
\rangle^{k+\frac{1}{2}-2\epsilon'}}  d\xi_1 d\xi_2 d\tau_1 
d\tau_2 \\
& \le & c \|v_1\|_{L^2_{xt}} \|v_2\|_{L^2_{xt}} \|\varphi\|_{L^2_{xt}} \, .   
\end{eqnarray*}
similarly as in Case 1. The regions where $|\sigma|$ or $|\sigma_1|$ are 
dominant can be handled similarly. \\
{\bf Case 3:} $|\xi_1| \sim |\xi_2|$ ($\Rightarrow \, |\xi| \le |\xi_1|+|\xi_2| 
\sim 2|\xi_2|$) \\
We use
$$ \frac{\langle \xi_1 \rangle^l}{\langle \xi \rangle^k \langle \xi_2 \rangle^l} 
\sim \frac{1}{\langle \xi \rangle^k}$$
and get in the $|\sigma_2|$ - dominant region (the other cases can be treated 
similarly again) by our assumption $k>0$ :
\begin{eqnarray*}
\lefteqn{\int \int_* 
\frac{|\tilde{v}_1(\xi_1,\tau_1)||\tilde{v}_2(-\xi_2,-\tau_2)| 
|\tilde{\varphi}(\xi,\tau)|}{ 
\langle \sigma_1 \rangle^{\frac{1}{2}+\epsilon'}
 \langle \sigma 
\rangle^{\frac{1}{2}+\epsilon'} \langle \xi_2 \rangle^{\frac{1}{2}-2\epsilon'} 
\langle \xi \rangle^k}  d\xi_1 d\xi_2 d\tau_1 
d\tau_2} \\
& \le & \int \int_* 
\frac{|\tilde{v}_1(\xi_1,\tau_1)||\tilde{v}_2(-\xi_2,-\tau_2)| 
|\tilde{\varphi}(\xi,\tau)|}{ \langle \sigma_1 \rangle^{\frac{1}{2}+\epsilon'} 
\langle \sigma \rangle^{\frac{1}{2}+\epsilon'} \langle \xi 
\rangle^{k+\frac{1}{2}-2\epsilon'}} d\xi_1 d\xi_2 d\tau_1 d\tau_2 \\
& \le & \left\| {\cal F}^{-1} \left( \frac{|\tilde{v}_1|}{\langle \sigma_1 
\rangle^ 
{\frac{1}{2}+\epsilon'}} \right) \right\|_{L^{\infty}_t L^2_x} \|v_2\|_{L^2_t 
L^2_x} \left\| {\cal 
F}^{-1} \left( \frac{|\tilde{\varphi}|}{\langle \xi \rangle^ 
{k+\frac{1}{2}-2\epsilon'}} \right) \right\|_{L^2_t L^{\infty}_x} 
   \\
& \le & c \|v_1\|_{L^2_{xt}} \|v_2\|_{L^2_{xt}} \|\varphi\|_{L^2_{xt}} \, .   
\end{eqnarray*}

{\bf B:} Consider now the (+,+) - or (--,--) - case, where one has to estimate
$$\int \int\limits_{* \atop \xi_1 \xi_2 > 0} 
\frac{|\tilde{v}_1(\xi_1,\tau_1)||\tilde{v}_2(-\xi_2,-\tau_2)| 
|\tilde{\varphi}(\xi,\tau)| \langle \xi_1 \rangle^l}{ 
\langle \tau_1 \pm |\xi_1| \rangle^{\frac{1}{2}+\epsilon'} \langle \tau_2 \mp 
|\xi_2| \rangle^{\frac{1}{2}-2\epsilon'} \langle \tau + |\xi| 
\rangle^{\frac{1}{2}+\epsilon'} \langle \xi \rangle^k \langle \xi_2 \rangle^l} 
d\xi_1 d\xi_2 
d\tau_1 d\tau_2 \, . $$
One has $|\xi| = |\xi_1| + |\xi_2|$ and with 
$$ \sigma_1 = \tau_1 \pm|\xi_1| \, , \, \sigma_2 = \tau_2 \mp|\xi_2| \, , \, 
\sigma = \tau + |\xi|$$
one checks again (\ref{*****}) .\\
If $|\sigma_2|$ is dominant and $|\xi_1| \ge |\xi_2|$ we have $|\xi| 
\sim|\xi_1|$ and
$$ \frac{\langle \xi_1 \rangle^l}{\langle \xi \rangle^k \langle \xi_2 \rangle^l} 
\sim \frac{1}{\langle \xi_1 \rangle^{k-l} \langle \xi_2 \rangle^l}$$
as well as $\langle \xi_2 \rangle \le c \langle \sigma_2 \rangle$ , so that by 
use of $k-l \ge 0$ , $k>0$ and $|\xi_1| \ge |\xi_2|$ we estimate
\begin{eqnarray*}
\lefteqn{\int \int_* 
\frac{|\tilde{v}_1(\xi_1,\tau_1)||\tilde{v}_2(-\xi_2,-\tau_2)| 
|\tilde{\varphi}(\xi,\tau)|}{ 
\langle \sigma_1 \rangle^{\frac{1}{2}+\epsilon'}
 \langle \sigma 
\rangle^{\frac{1}{2}+\epsilon'} \langle \xi_1 \rangle^{k-l} \langle \xi_2 
\rangle^{l+\frac{1}{2}-2\epsilon'}}  d\xi_1 d\xi_2 d\tau_1 
d\tau_2} \\
& \le & \int \int_* 
\frac{|\tilde{v}_1(\xi_1,\tau_1)||\tilde{v}_2(-\xi_2,-\tau_2)| 
|\tilde{\varphi}(\xi,\tau)|}{ \langle \sigma_1 \rangle^{\frac{1}{2}+\epsilon'} 
\langle \sigma \rangle^{\frac{1}{2}+\epsilon'} \langle \xi_2 
\rangle^{k+\frac{1}{2}-2\epsilon'}} d\xi_1 d\xi_2 d\tau_1 d\tau_2 \\
& \le & \left\| {\cal F}^{-1} \left( \frac{|\tilde{v}_1|}{\langle \sigma_1 
\rangle^ 
{\frac{1}{2}+\epsilon'}} \right) \right\|_{L^{\infty}_t L^2_x}  \left\| {\cal 
F}^{-1} \left( \frac{|\tilde{v}_2(-\xi_2,-\tau_2)|}{\langle \xi_2 \rangle^ 
{k+\frac{1}{2}-2\epsilon'}} \right) \right\|_{L^2_t L^{\infty}_x} 
\|\varphi\|_{L^2_t L^2_x} 
   \\
& \le & c \|v_1\|_{L^2_{xt}} \|v_2\|_{L^2_{xt}} \|\varphi\|_{L^2_{xt}} \, .   
\end{eqnarray*}
If $|\sigma_2|$ is dominant and $|\xi_2| \ge |\xi_1|$, we have $|\xi| \sim 
|\xi_2|$ and, using $k+l>0$ :
$$ \frac{ \langle \xi_1 \rangle^l}{\langle \xi \rangle^k \langle \xi_2 
\rangle^l} \sim \frac{ \langle \xi_1 \rangle^l}{ \langle \xi_2 \rangle^{k+l}} 
\le  \frac{ \langle \xi_1 \rangle^l}{ \langle \xi_1\rangle^{k+l}} = 
\frac{1}{\langle \xi_1 \rangle^k} $$
and also $\langle \xi_1 \rangle \le c \langle \sigma_2 \rangle \, . $
Thus, similarly as before we get for $k>0$ :
\begin{eqnarray*}
& & \int \int_* 
\frac{|\tilde{v}_1(\xi_1,\tau_1)||\tilde{v}_2(-\xi_2,-\tau_2)| 
|\tilde{\varphi}(\xi,\tau)|}{ 
\langle \sigma_1 \rangle^{\frac{1}{2}+\epsilon'}
 \langle \sigma 
\rangle^{\frac{1}{2}+\epsilon'} \langle \xi_1 
\rangle^{k+\frac{1}{2}-2\epsilon'}}  d\xi_1 d\xi_2 d\tau_1 
d\tau_2 \\
& &  \le  c \|v_1\|_{L^2_{xt}} \|v_2\|_{L^2_{xt}} \|\varphi\|_{L^2_{xt}} \, .
\end{eqnarray*}
The other regions are treated similarly.\\[1em]
{\bf Remark.} The modified estimate (\ref{***}) with $ 
Y^{-k,-\frac{1}{2}-\epsilon'}_+$ replaced by  $ Y^{-k,-\frac{1}{2}-\epsilon'}_-$ 
is proven in the same way. See also the remark to the previous lemma.

We summarize our results in the following
\begin{theorem}
\label{Theorem 1}
Assume $ l < \frac{1}{4}$ , $ k > 0 $ , $2l+k <1$ , $ l+k \le 1$ and $ k \ge |l| 
$ . The Cauchy problem for the Dirac -- Klein -- Gordon equations 
(\ref{0.1}),(\ref{0.2}),(\ref{0.3}) with data
$$\psi_0 \in H^{-l}({\bf R}) \,  , \, \phi_0 \in H^k({\bf R}) \, , \, \phi_1 \in 
H^{k-1}({\bf R}) $$
has a unique local solution
$$\psi = \psi_+ + \psi_- \quad
{\mbox with} \quad \psi_{\pm} \in X^{-l,\frac{1}{2}+\epsilon'}_{\pm}[0,T] \, ,$$ 
and
$$\phi = \frac{1}{2}(\phi_+ + \phi_-)\, , 
\, \phi_t = \frac{1}{2i} A^{\frac{1}{2}}(\phi_+ - \phi_-) \quad
{\mbox with} \quad \phi_{\pm} \in Y^{k,\frac{1}{2}+\epsilon'}_{\pm}[0,T] \, ,$$  
where $A=-\frac{\partial^2}{\partial x^2} +1$ . Here 
$T=T(\|\psi_0\|_{H^{-l}},\|\phi_0\|_{H^k},\|\phi_1\|_{H^{k-1}})$ and $\epsilon' 
> 0$ is sufficiently small. This solution satisfies 
$$\psi \in C^0([0,T],H^{-l}({\bf R})) \, , \, \phi \in C^0([0,T],H^k({\bf R})) 
\, , \, \phi_t \in C^0([0,T],H^{k-1}({\bf R})) \, .$$ 
\end{theorem} 
\section{Global existence}
The following global existence result is an easy consequence of the local 
results and conservation of charge.
\begin{theorem}
\label{Theorem 2}
Assume $\psi_0 \in L^2({\bf R})$ , $ \phi_0 \in H^k({\bf R}) $ , $ \phi_1 \in 
H^{k-1}({\bf R})$ , where $ 0 < k < \frac{1}{2} $ . Then the solution of Theorem 
\ref{Theorem 1} exists globally in $t$.
\end{theorem}
{\bf Proof:} We only need an a-priori-bound for $\|\psi(t)\|_{L^2}$ and 
$\|\phi(t)\|_{H^k} + \|\phi_t(t)\|_{H^{k-1}}.$ Charge conservation gives the 
$L^2$ - bound of $\psi(t)$ and $\phi(t)$ fulfills the integral equation
\begin{eqnarray*}
\phi(t) & = & \cos(A^{\frac{1}{2}}t) \phi_0 + A^{-\frac{1}{2}} 
\sin(A^{\frac{1}{2}}t)\phi_1 + \int_0^t A^{-\frac{1}{2}} 
\sin[A^{\frac{1}{2}}(t-s)] \langle \beta \psi(s),\psi(s) \rangle ds \\
& & \hspace{6cm} + c_0  \int_0^t A^{-\frac{1}{2}} \sin[A^{\frac{1}{2}}(t-s)] 
\phi(s) ds \, ,
\end{eqnarray*}
where $c_0 = 1-m^2$ . Thus
\begin{eqnarray*}
\lefteqn{ \|\phi(t)\|_{H^k} + \|\phi_t(t)\|_{H^{k-1}} }\\
& & \le c(\|\phi_0\|_{H^k} + \|\phi_1\|_{H^{k-1}} + \int_0^t (\| \langle \beta 
\psi(s),\psi(s) \rangle \|_{H^{k-1}} + \|\phi(s)\|_{H^{k-1}}) \, ds) \, .
\end{eqnarray*}
Using the estimate
$$ \| \langle \beta \psi , \psi \rangle \|_{H^{k-1}} \le c \|\psi \|_{L^2}^2 
\quad {\mbox for} \quad k < \frac{1}{2} \, , $$
which follows from (cf. \cite{BG})
\begin{eqnarray*}
\|uv\|_{H^{k-1}}^2 & \le & \int \left| \int \tilde{u}(\eta) \tilde{v}(\xi - 
\eta) \, d\eta \right|^2 \langle \xi \rangle^{2(k-1)} \, d\xi \\
& \le & \int ( \int |\tilde{u}(\eta)|^2 d\eta)(\int|\tilde{v}(\xi - \eta)|^2 
d\eta) \langle \xi \rangle^{2(k-1)} d\xi \\
& \le & \|u\|_{L^2}^2 \|v\|_{L^2}^2 \int \langle \xi \rangle^{2(k-1)} d\xi \le c 
\|u\|_{L^2}^2 \|v\|_{L^2}^2 \, ,
\end{eqnarray*}
we arrive at
$$ \|\phi(t)\|_{H^k} + \|\phi_t(t)\|_{H^{k-1}} \le c(\|\phi_0\|_{H^k} + 
\|\phi_1\|_{H^{k-1}} + t \|\psi_0\|_{L^2}^2 + \int_0^t \|\phi(s)\|_{H^k} \, ds) 
\, , $$
so that Gronwall's lemma gives the desired a-priori-bound.

\end{document}